# The Wiener Index and the Hosoya Polynomial of the Jahangir Graphs

Shaohui Wang[1], Mohammad Reza Farahani[2] [*], M. R. Rajesh Kanna[3],
Muhammad Kamran Jamil[4], R. Pradeep Kumar[5]

[1]Department of Mathematics, University of Mississippi, University, MS, USA

[1]Department of Mathematics and Computer Science, Adelphi University, Garden City, NY, USA

[2]Department of Applied Mathematics, Iran University of Science and Technology (IUST) Narmak, Tehran, Iran

[3]Department of Mathematics, Maharani's Science College for Women, Mysore, India

[4]Department of Mathematics, Riphah Institute of Computing and Applied Sciences (RICAS), Riphah International University, Lahore, Pakistan

[5]Department of Mathematics, the National Institute of Engineering, Mysuru, India

**Email address:**
shaohuiwang@yahoo.com (S. Wang), mrfarahani88@gmail.com (M. R. Farahani), mr.rajeshkanna@gmail.com (M. R. R. Kanna), m.kamran.sms@gmail.com (M. K. Jamil), pradeepr.mysore@gmail.com (P. P. Kumar).

[*]Corresponding author

**Abstract:** Let G be a simple connected graph having vertex set V and edge set E. The vertex-set and edge-set of G denoted by V(G) and E(G), respectively. The length of the smallest path between vertices u,v ∈ V(G) is called the distance, d(u,v), between the vertices u,v. Mathematical chemistry is the area of research engaged in new application of mathematics in chemistry. In mathematics chemistry, we have many topological indices for any molecular graph, that they are invariant on the graph automorphism. In this research paper, we computing the Wiener index and the Hosoya polynomial of the Jahangir graphs $J_{5,m}$ for all integer number m≥3. The Wiener index is the sum of distances between all pairs of vertices of G as W(G)= $\sum_{\{u,v\}\subset V(G)} d(u,v)$. And the Hosoya polynomial of G is H(G,x)= $\sum_{\{u,v\}\subset V(G)} x^{d(u,v)}$, where d(u,v) denotes the distance between vertices u and v.

## 1. Introduction

Let G be a connected graph. The vertex-set and edge-set of G denoted by *V(G)* and *E(G)*, respectively. The degree of a vertex $v \in V(G)$ is the number of vertices joining to *v* and denoted by $d_v$

Topological indices of a simple graph are numerical descriptors that are derived from graph of chemical compounds. Such indices based on the distances in graph are widely used for establishing relationships between the structure of molecular graphs and Nanotubes and their physicochemical properties. Usage of topological indices in biology and chemistry began in 1947 when chemist *Harold Wiener* [1] introduced the *Wiener index* to demonstrate correlations between physicochemical properties of organic compounds and the index of their molecular graphs. Wiener originally defined his index on trees and studied its use for correlations of physico-chemical properties of alkanes, alcohols, amines and their analogous compounds [2] as:

$$W(G) = \tfrac{1}{2} \sum_{v \in V(G)} \sum_{u \in V(G)} d(u,u) \qquad (1)$$

where *d(u,v)* denotes the distance between vertices *u* and *v*.

The Hosoya polynomial of a graph is a generating function about distance distributing, introduced by *Haruo Hosoya* in 1988 and for a connected graph *G* is defined as [3]:

$$H(G,x) = \tfrac{1}{2} \sum_{v \in V(G)} \sum_{u \in V(G)} x^{d(v,u)} \qquad (2)$$

In a series of papers, the Wiener index and the Hosoya polynomial of some molecular graphs and Nanotubes are computed. For more details about the Wiener index and the Hosoya polynomial, please see the paper series [4-22] and the references therein.

In this research paper, we present some properties of the Wiener index and the Hosoya polynomial and we introduce a closed formula of this and the correspondent polynomial of the Jahangir *graphs $J_{5,m}$* for all integer number *m≥3*.

## 2. Materials and Methods

The In the this section, we present some studies about a semi-regular and connected graphs that named "Jahangir Graphs $J_{n,m}$" $\forall n \geq 2$, $\forall m \geq 3$. And we introduce a method to compute the Wiener index and the Hosoya polynomial of the Jahangir graphs $J_{5,m}$ in continue.

What is a the Jahangir graphs $J_{n,m}$?

Definition 2.1. [23, 24] The Jahangir graphs $J_{n,m}$ is a graph on nm+1 vertices and m(n+1) edges $\forall n \geq 2$ & $\forall m \geq 3$; i.e., a graph consisting of a cycle $C_{nm}$ with one additional vertex which is adjacent to m vertices of $C_{nm}$ at distance n to each other on $C_{nm}$.

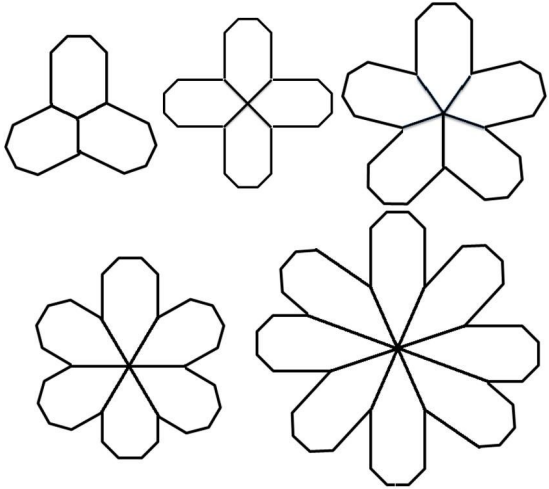

**Figure 1.** Some examples of Jahangir graphs $J_{5,3}, J_{5,4}, J_{5,5}, J_{5,6}$ and $J_{5,8}$.

In particular, $\forall m \geq 3$ the Jahangir graphs $J_{5,m}$ is a graph consisting of a cycle $C_{5m}$ with one additional vertex (the Center vertex c) which is adjacent to m vertices of $C_{nm}$ at distance 5 to each other on $C_{5m}$. Some first example of this graphs are shown in Figure 1. For more details about the Jahangir graphs $J_{n,m}$ reader can see the paper series [23-33,34-36].

## 3. Main Results and Discussion

In this paper, we computed a closed formula of the Wiener index and the Hosoya polynomial of the Jahangir *graphs* $J_{5,m}$ for all integer number $m \geq 3$.

Theorem 2.1. The Hosoya polynomial and the Wiener index for the Jahangir graphs $J_{5,m}$ for all integer number $m \geq 3$ are equal to

$$H(J_{5,m},x) = 6mx^1 + \tfrac{1}{2}(m^2+13m)x^2 + (2m^2+5m)x^3 + (4m^2-4m)x^4$$
$$+ (4m^2-6m)x^5 + (2m^2-5m)x^6 \quad (3)$$

$$W(J_{5,m}) = 55m^2 - 42m. \quad (4)$$

*Proof.* Let $J_{5,m}$ be Jahangir graphs $\forall m \geq 3$ with 5m+1 *vertices and 6m* edges. From Definition 2.1 and Figure 1, one can see that there are *4m* vertices of the Jahangir graph $J_{5,m}$ with degree 2 and *m* vertices of $J_{5,m}$ with degree 3 and Center vertex *c* has degree *m* and we have three partitions of the vertex set $V(J_{5,m})$ as follow

$$V_2 = \{v \in V(J_{5,m}) | d_v = 2\} \rightarrow |V_2| = 4m \quad (5)$$

$$V_3 = \{v \in V(J_{5,m}) | d_v = 3\} \rightarrow |V_3| = m \quad (6)$$

$$V_m = \{c \in V(J_{5,m}) | d_c = m\} \rightarrow |V_m| = 1 \quad (7)$$

And alternatively, $V(J_{5,m}) = V_2 \cup V_3 \cup V_m$ and $V_2 \cap V_3 \cap V_m = \emptyset$ and we know

$$|E(G)| = \tfrac{1}{2} \sum_{k=\delta}^{\Delta} |V_k| \times k \quad (8)$$

where $\delta$ and $\Delta$ are the minimum and maximum of $d_v$ for all $v \in V(G)$, respectively, thus

$$|E(J_{5,m})| = \tfrac{1}{2}[2 \times |V_2| + 3 \times |V_3| + 3 \times |V_m|]. \quad (9)$$

Now, for compute the Hosoya polynomial and the Wiener index of $J_{5,m}$, we denote the number of unordered pairs of vertices u and v of a graph G as distance d(u,v)=k by d(G,k) for all integer number k up to d(G) (where d(G) denote the topological diameter and is the longest distance between vertices of a graph G). Thus, we redefine this mention topological polynomial and index of G as follow:

$$H(G,x) = \sum_{k=1}^{d(G)} d(G,k) x^k \quad (10)$$

$$W(G) = \sum_{k=1}^{d(G)} d(G,k) \times k \quad (11)$$

From the Definition 2.1 and the structure of Jahangir graphs $J_{5,m}$ in Figure 1, we have following computations for $d(J_{5,m}, k)$ ($\forall k: 1 \leq k \leq d(J_{5,m})$):

$$d(J_{5,m},1) = |E(J_{5,m})| = 6m \text{ by definitions of } d(G,k). \quad (12)$$

$$d(J_{5,m},2) = \tfrac{1}{2} m(m+13) \quad (13)$$

Since, there are two 2-edges paths between Center vertex $c \in V(J_{5,m})$ and other vertices of vertex set $V_2 \subset V(J_{5,m})$). $\tfrac{1}{2}|V_3|(m-1)$ 2-edges paths between all vertices of $u,v \in V_3 \subset V(J_{5,m})$), and $2|V_3| + |V_3|$ 2-edges paths start from vertices of $V_2$ until vertices of $V_3$ and $V_2 \subset V(J_{5,m})$. Thus the second sentence of the Hosoya polynomial $H(J_{5,m},x)$ of $J_{5,m}$ is equal to $[2|V_3| + \tfrac{1}{2}|V_3|(m-1) + 3|V_3| + 2|V_3|]x^2 = \tfrac{1}{2}(m^2 + 13m)x^2$.

$$d(J_{5,m},3) = (2m^2 + 5m) \quad (14)$$

Because, there are $2|V_3|$ 3-edges paths between $c \in V_m$ and vertices of $V_2$, but there are not any 3-edges paths between $c$ and vertices of $V_3$. Also, there are $\tfrac{1}{2} \times 2|V_3| + 2|V_3| = 3m$ 3-edges paths between all vertices of $u,v \in V_2 \subset V(J_{5,m})$). From Figure 1, one can see that $2|V_3| + |V_3|(2(m-1))$ 3-edges paths started from vertices of $V_3$ until vertices of $V_2 \subset V(J_{5,m})$. Thus, the third sentence of $H(J_{5,m},x)$ is $\tfrac{1}{2}[2|V_3| + 0 + 3|V_3| + 2m|V_3|]x^3 = m(2m+5)x^3$.

$$d(J_{5,m},4) = 4m^2 - 2m \quad (15)$$

Since from Figure 1, there are $2|V_3|$ 3-edges paths between $c \in V_m$ and vertices of $V_2$. And obviously there are not any 4-

edges paths between vertices of $V_3$ and $c$. And there are $|V_3|(2|V_3|-4)$ 4-edges paths between vertices of $V_3$ and $V_2$. Also, there are $d(C_m)+½×2|V_3|(2|V_3|-3)$ 4-edges paths between all vertices of $u,v \in V_2 \subset V(J_{5,m})$). Therefore the coefficient of the fourth sentence of $H(J_{5,m},x)$ is $2|V_3|+0+|V_3|(2|V_3|-4)+2m(m-1)=4m(m-1)=4m^2-4m$.

$$d(J_{5,m},5)=(4m^2-6m) \quad (16)$$

From the Definition 1 and the structure of $J_{5,m}$, one can see that there are not any 5-edges paths between vertices of $V_3$ and $c$ and vertices of $V_2$ and between vertices of $u \in V_2$ and $v \in V_3$. Thus for the $5^{th}$ sentence of the Hosoya polynomial of $J_{5,m}$, we have $2d(C_m)+½×|V_2|(2|V_3|-4)=2m(2m-3)$ 5-edges paths between all vertices of $u,v \in V_2$.

$$\text{For } d(J_{5,m})=6; \ d(J_{5,m},6)=(2m^2-5m) \quad (17)$$

Because, there isn't any 6-edges paths started from members of $C$ and $V_3$ to all other vertices of Jahangir graphs $J_{5,m}$ in Figure 1. But there are $½×|V_2|(2|V_2|-5)=m(2m-5)$ 6-edges paths between all vertices $V_2$. Thus, the 6th and last sentence of $H(J_{5,m},x)$ is $m(2m-5)x^6$.

Now, Equations 12, 13,…,17 imply that the Hosoya polynomial of the Jahangir graph $J_{5,m}$ is equal to:

$$H(J_{5,m},x)=6mx^1+½(m^2+13m)x^2+(2m^2+5m)x^3 \\ +(4m^2-4m)x^4+(4m^2-6m)x^5+(2m^2-5m)x^6 \quad (18)$$

Also, from the definition of Wiener index and the Hosoya Polynomial of a graph $G$ and the first derivative of Hosoya polynomial (evaluated at $x=1$ in Equations 10, 11), we cam compute the Wiener index of the Jahangir graph $J_{5,m}$ as follow:

$$W(J_{5,m})=H(J_{5,m},x)'|_{x=1}$$

$$=\sum_{k=1}^{d(J_{5,m})} d(J_{5,m},k) \times k$$

$$=[6mx^1+½(m^2+13m)x^2+(2m^2+5m)x^3+(4m^2-4m)x^4+(4m^2-6m)x^5+(2m^2-5m)x^6]'|_{x=1}$$

$$=[6m×1+½m(m+13)×2+m(2m+5)×3+4m(m-1)×4+m(4m-6)×5+m(2m-5)×6]$$

$$=55m^2-42m. \quad (19)$$

Here the proof of theorem is completed.

Example 2.1. The Hosoya polynomial and the Wiener index for $J_{5,6}$ are equal to:

$$H(J_{5,6},x)=36x+57x^2+102x^3+120x^4+108x^5+42x^6 \quad (20)$$

$$W(J_{5,6})=55(6)^2-42(6)=1728 \quad (21)$$

## 4. Conclusion

In this paper, we obtained the Hosoya polynomial and the Wiener index of graph "the Jahangir *graphs* $J_{5,m}$ for all integer number $m \geq 3$" for the first time as:

$$H(J_{5,m},x)=6mx^1+½(m^2+13m)x^2+(2m^2+5m)x^3+(4m^2-4m)x^4 \\ +(4m^2-6m)x^5+(2m^2-5m)x^6$$

$$W(J_{5,m})=55m^2-42m.$$

## Acknowledgements

The author is thankful to Professor Emeric Deutsch from Department of Mathematics of Polytechnic University (Brooklyn, NY 11201, USA) for his precious support and suggestions.